\documentclass[a4paper]{amsart}
\usepackage{a4wide}
\usepackage{amsmath,amssymb,mathrsfs}
\usepackage[UKenglish] {babel}

\def\a{\alpha} \def\b{\beta} \def\d{\delta} \def\e{\epsilon} \def\f{\varphi}
\def\l{\lambda}  \def\R{\mathbb{R}} \def\C{\mathbb{C}}
\def\H{\mathbb{H}} 
\def\O{\mathbb{O}} \def\I{\mathbb{I}}
 \def\minus{\smallsetminus}
\def\({\left(} \def\){\right)} 
\def\<{\langle} \def\>{\rangle}

\def\inv{^{-1}}

\renewcommand\ge{\geqslant}
\renewcommand\le{\leqslant}
\newcommand\ie{i.e.}

\newcommand\eg{e.g.}

\newcommand\forget[1]{}

\newcommand{\ct}{\mathrm{C}_2}

\DeclareMathOperator{\spann}{span}
\DeclareMathOperator{\rt}{R}
\DeclareMathOperator{\lt}{L}
\DeclareMathOperator{\sign}{sign}

\DeclareMathOperator{\Aut}{Aut}
\DeclareMathOperator{\Iso}{Iso}
\DeclareMathOperator{\og}{O}
\DeclareMathOperator{\so}{SO}
\DeclareMathOperator{\Pds}{Pds}
\DeclareMathOperator{\spds}{SPds}
\DeclareMathOperator{\SL}{SL}
\DeclareMathOperator{\gl}{GL}

\DeclareMathOperator{\pgo}{PGO}
\DeclareMathOperator{\go}{GO}

\DeclareMathOperator{\pgl}{PGL}

\DeclareMathOperator{\g}{G}
\DeclareMathOperator{\pg}{PG}

\newcommand{\chr}{\mathop{\rm char}\nolimits}

\newenvironment{smatrix}{\left(\begin{smallmatrix}}{\end{smallmatrix}\right)}

\newtheorem{thm}{Theorem}
\newtheorem{lma}[thm]{Lemma}
\newtheorem{prop}[thm]{Proposition}
\newtheorem{cor}[thm]{Corollary}
\theoremstyle{definition}

\theoremstyle{remark}
\newtheorem{rmk}[thm]{Remark}

\title{Isotopes of Hurwitz algebras}
\author{Erik Darp\"o}
\address{Graduate School of Mathematics, Nagoya University, Furo-cho, Chikusa-ku, Nagoya, Japan}

\begin{document}
\selectlanguage{UKenglish}
\date{}

\begin{abstract}
We study the class of all algebras that are isotopic to a Hurwitz algebra.
In the general case, isomorphism classes of such algebras are shown to correspond to
orbits of a certain group action. Our main focus is the isotopes of Hamilton's
quaternions, for which a more detailed, and geometrically intuitive description of the
isomorphism classes is given. 
As an application, we also demonstrate how some results concerning the classification of
finite-dimensional composition algebras can be deduced from our general results. 
\end{abstract}
\thanks{\emph{MSC~2010:} 17A60 (17A35, 17A75)}

\maketitle 

\section{Introduction} \label{introduction}
Let $k$ be a field, and $V$ a vector space over $k$.
We shall say that a quadratic form $q:V\to k$ is \emph{non-degenerate} if the associated
bilinear form $\<x,y\>=q(x+y)-q(x)-q(y)$ is non-degenerate (\ie, $\<x,V\>=0$ only if
$x=0$).
A \emph{composition algebra} is a non-zero (not necessarily associative) algebra $A$ over
a field $k$, equipped with a non-degenerate quadratic form $n:A\to k$ such that
$n(ab)=n(a)n(b)$ for all $a,b\in A$. The form $n$ is usually called the \emph{norm} of
$A$. A composition algebra possessesing an identity element is called a 
\emph{Hurwitz algebra}.
Every Hurwitz algebra has dimension one, two, four or eight (thus, in particular, it is
finite dimensional), and can be constructed via an iterative method known as the
\emph{Cayley-Dickson process}.\footnote{Some authors have used a weaker notion of
  non-degeneracy for the norm $n$, requiring that $n(x+y)=n(y)$ for all $y$ implies
  $x=0$. This definition gives rise to additional unital composition algebras over fields
  $k$ of characteristic two, in form of purely inseparable field extensions of $k$
  \cite{kaplansky53}. If $\chr k\ne2$, the two definitions are equivalent.}
The facts about Hurwitz algebras referred in this section are described in detail in
\cite{involutions}, Chapter~VIII. 

Two algebras $A$ and $B$ over a field $k$ are said to be \emph{isotopic} if there exist
invertible linear maps $\a,\b,\gamma:A\to B$ such that $\gamma(xy)=\a(x)\b(y)$ for all
$x,y\in A$. Clearly, isotopy is an equivalence relation amongst $k$-algebras. If $A$ and $B$
are isotopic then there exist $\a,\b\in\gl(A)$ such that the algebra $(A,\circ)$, with
multiplication $x\circ y=\a(x)\b(y)$, is isomorphic to $B$.
The algebra $(A,\circ)$ is called the \emph{principal isotope} of $A$
determined by $\a$ and $\b$, and is denoted by $A_{\a,\b}$. 

Several important classes of non-associative algebras can be constructed by isotopy from
the Hurwitz algebras. Examples include:
\begin{enumerate}
\item All finite-dimensional composition algebras \cite[p.~957]{kaplansky53}. This
  includes in particular all 
  finite-dimensional absolute valued algebras, which are precisely the finite-dimensional
  composition algebras over $\R$ whose norm is anisotropic (\ie, $n(x)=0$ only if $x=0$). 
  However, there exist composition algebras without identity element that are
  infinite dimensional, and thus not isotopic to any Hurwitz algebra; see \eg{}
  \cite{cuenca92,eldmyung93}.
\item All division algebras of dimension two over a field of characteristic different from
  two \cite{petersson00}.
\item All finite-dimensional division algebras $A$ that are not isotopic to an associative
  algebra, and satisfy the following property: 
  for all non-zero $a\in A$ there exists an element $b\in A$ such that $b(ax)=x$ for all
  $x\in A$.
  \cite{invqgi}
\end{enumerate}

The purpose of the present article is to describe, in a uniform way, all isotopes of
Hurwitz algebras. 
Generalising ideas that have earlier been used in more specialised  
situations (for example in \cite{ckmmrr11,petersson00}), we give a general description of 
all algebras isotopic to a Hurwitz algebra, encompassing also the case of the ground field
$k$ having characteristic two. 
Our principal result is a comprehensive geometric description of all 
isotopes of Hamilton's quaternion algebra $\H$ -- a class of real division algebras that
has not been studied before.

Isotopes of Hurwitz algebras of dimension two is a special case of Petersson's
\cite{petersson00} description of all $2$-dimensional algebras, which is also based on
the concept of isotopy. There is a large number of articles on finite-dimensional
composition algebras and absolute valued algebras (some examples are
\cite{alsaody15,ckmmrr11,eldmyung93,kaplansky53,ramirez99}), most of which use the concept of
isotopy as a central tool.

This article is organised as follows.
In Section~\ref{general}, a general description is given of the category of isotopes of a
Hurwitz algebra $A$. A more elaborate study of the case where $(A,n)$ is a Euclidean
space is given in Section~\ref{euclidean}, bringing about the announced description of
isotopes of $\H$. Finally, Section~\ref{calg} treats composition algebras,
showing how a description of these can be deduced from the results in
Section~\ref{general}.

From here on, let $k$ denote a field. All algebras are, unless otherwise stated, assumed
to be finite dimensional over $k$. Every element $a$ of an algebra $A$ determines linear
endomorphisms $\lt_a$ and $\rt_a$ of $A$, defined by $\lt_a(x)=ax$ and $\rt_a(x)=xa$
respectively. For a principal isotope $A_{\a,\b}=(A,\circ)$ of $A$, we use the notation
$\lt^\circ_a(x)=a\circ x$ and $\rt^\circ_a(x)=x\circ a$.
An algebra $A$ is said to be a \emph{division algebra} if $ A\ne0$ and
$\lt_a$ and $\rt_a$ are bijective for all non-zero $a\in A$. Moreover, $A$ is
\emph{alternative} if the identities $x^2y=x(xy)$ and $xy^2=(xy)y$ hold for all $x,y\in A$.

All Hurwitz algebras are alternative. 
Two-dimensional Hurwitz algebras are quadratic \'etale algebras, \ie, either 
separable field extensions of $k$ or isomorphic to $k\times k$. 
The Hurwitz algebras of dimension four are all \emph{quaternion algebras}, that is,
all four-dimensional central simple associative algebras. Eight-dimensional Hurwitz
algebras are precisely the central simple alternative algebras that are not associative
\cite{zorn31} (these are called \emph{octonion algebras}).
A Hurwitz algebra $A$ is commutative if and only if $\dim A\le2$, and associative if and
only if $\dim A\le4$. 

Any element $x$ in a Hurwitz algebra $A=(A,n)$ satisfies $x^2=\<x,1\>x-n(x)1$.
Hence, the norm $n$ is uniquely determined by the algebra structure of $A$, and every
algebra morphism of Hurwitz algebras that respects the identity element also preserves
the norm.  
Every Hurwitz algebra of dimension at least $2$ has a unique non-trivial
involution\footnote{By an involution is meant a self-inverse isomorphism $A\to A^{\rm op}$.} 
$\kappa:A\to A,\:x\mapsto\bar{x}$ satisfying $x+\bar{x}\in k1$ and
$x\bar{x}=\bar{x}x=n(x)1$ for all $x\in A$. 
Moreover, two Hurwitz algebras are isomorphic if and only if their respective norms are
equivalent (this was first proved in \cite{jacobson58} in characteristic different from
two). The quadratic forms occurring as norms of Hurwitz algebras are precisely the
$m$-fold Pfister forms over $k$, for $m\in\{0,1,2,3\}$.
If $A$ is a Hurwitz algebra and $a\in A$, then $\lt_a$ and $\rt_a$ are invertible
if and only if $n(a)\ne0$. This is also equivalent to the existence of an inverse $a\inv$
of $a$ in $A$: since $a\bar{a}=\bar{a}a=n(a)1$, we have $a\inv=n(a)\inv\bar{a}$ if
$n(a)\ne0$. 
Moreover, $\lt_a\inv=\lt_{a\inv}$ and $\rt_a\inv=\rt_{a\inv}$ in this case.
The invertible elements of any alternative algebra $A$ form a so-calles 
\emph{Moufang loop} under multiplication (the concept was introduced by Moufang in
\cite{moufang35} under the name \emph{quasi-group}), denoted by $A^\ast$.

For any algebra $A$, the \emph{nucleus} is defined as
$N(A)=\{a\in A\mid (xy)z=x(yz) \mbox{ for } a\in\{x,y,z\}\}$. The nucleus is an associative
subalgebra of $A$. If $A$ is a Hurwitz algebra, then $a\in N(A)$ if and only if
$(xa)y=x(ay)$ for all $x,y\in A$. If $\dim A\le4$ then $A$ is associative and thus
$N(A)=A$; the nucleus of an eight-dimensional Hurwitz algebra is $k1$.

A \emph{similitude} of a non-zero quadratic space $V=(V,q)$ is an invertible linear map
$\f:V\to V$ such that $q(\f(x))=\mu(\f)q(x)$ for all $x\in V$, where $\mu(\f)\in k$ is a
scalar independent of $x$. The element $\mu(\f)$ is called the \emph{multiplier} of $\f$.
If $l=\dim V$ is even, then $\det(\f)=\pm\mu(\f)^{l/2}$ \cite[12A]{involutions}.
If $\chr k\ne2$, a similitude $\f$ satisfying $\det(\f)=\mu(\f)^{l/2}$ are said to be
\emph{proper}. In the characteristic two case, a similitude $\f:V\to V$ is defined to be
proper if its Dickson invariant (see \cite[12.12]{involutions}) is zero.
The group of all similitudes of $V$ is denoted by $\go(V,q)$, or $\go(V)$ for short. The
proper similitudes form a normal subgroup $\go^+(V)\subset\go(V)$ of index two. Elements
in $\go(V)\minus\go^+(V)$ are called \emph{improper} similitudes.
The map $\mu:\go(V)\to k^\ast,\f\mapsto\mu(\f)$ is a group homomorphism, the kernel of
which is the orthogonal group $\og(V)$.
We write $\og^+(V)=\og(V)\cap\go^+(V)$, or $\so(V)=\og^+(V)$ in case $k=\R$ and $(V,q)$ is
a Euclidean space.

For any Hurwitz algebra $A=(A,n)$, define 
$\g(A)$ to be the set of all $\f\in\go(A)$ for which there exist $\f_1,\f_2\in\go(A)$ such
that $\f(xy)=\f_1(x)\f_2(x)$ for all $x,y\in A$.
Now, by what is known as the principle of \emph{triality} for Hurwitz algebras
\cite[3.2]{sv00},
$$\g(A)= \begin{cases}
  \go^+(A) & \mbox{if } \dim A\ge4, \\ \go(A) & \mbox{if } \dim A\le2.
\end{cases}$$
We call $\f_1$ and $\f_2$ \emph{triality components} of $\f$.
It is not difficult to see that any other pair of triality components of
$\f$ is of the form $(\rt_w\inv\f_1,\lt_w\f_2)$ for some $w\in N(A)^\ast$, and that
$\f_1=\rt_{\f_2(1)}\inv\f$ and $\f_2=\lt_{\f_1(1)}\inv\f$.
If $A$ is associative, then
$$\f(xy)=\f_1(x)\f_2(y)=(\f(x)\f_2(1)\inv)(\f_1(1)\inv\f(y))=\f(x)\f(1)\inv\f(y)\,.$$
Given a Hurwitz algebra $A$, triality components $\f_1,\f_2$ of any $\f\in\g(A)$ are
again elements of $\g(A)$.
Moreover, $(\f_1)\inv$ and $(\f_2)\inv$ are a pair of triality components of $\f\inv$.
We write $\pg(A)=\g(A)\slash(k^\ast\I)$. 

A \emph{groupoid} is a category in which every morphism is an isomorphism. 
Every action of a group $G$ on a set $X$ gives rise to a groupoid with object class $X$,
and morphisms $x\to y$ being the set of group elements $g\in G$ satisfying $g\cdot x=y$. 
We call this the \emph{groupoid of the $G$-action on $X$}.
Given a vector space $V$ over $k$, and a quadratic form $q:V\to k$, set
$\pgl(V)=\gl(V)/(k^\ast\I)$, $\pgo(V,q)=\go(V,q)/(k^\ast\I)$ and
$\pgo^+(V,q)=\go^+(V,q)/(k^\ast\I)$.
The orbit/coset (of some group action/normal subgroup) represented by an element in a
set/group will be denoted using square brackets around the element in question. For
example, if $\a\in\gl(A)$, the element in $\pgl(A)$ represented by $\a$ is written as
$[\a]$.
If $V$ is a Euclidean space then $\Pds(V)$ denotes the set of positive definite symmetric
endomorphisms of $V$, and $\spds(V)=\Pds(V)\cap\SL(A)$. The set of isomorphisms from an object $A$ to an object $B$ in a
category $\mathscr{C}$ is denoted by $\Iso(A,B)=\Iso_{\mathscr{C}}(A,B)$.
Throughout, $\ct$ denotes the cyclic group of order two, generated by the canonical
involution in a Hurwitz algebra of dimension at least two: $\ct=\<\kappa\>=\{\I,\kappa\}$.

\section{General description} \label{general}

The following lemma is essentially due to Albert \cite[Theorem~7]{albert42a}, for general
(not necessarily associative) algebras. The formulation has been sharpened slightly and
adapted to our situtation.

\begin{lma} \label{identity}
    Let $A$ be a Hurwitz algebra, and $\a,\b\in\gl(A)$. The isotope $A_{\a,\b}$ is unital
    if and only if $\a=\rt_a$, $\b=\lt_b$ for some $a,b\in A^\ast$.
    The identity element in $A_{\rt_a,\lt_b}$ is $(ab)\inv$.
  \end{lma}

\begin{proof}
The isotope $A_{\a,\b}=(A,\circ)$ is unital if and only if there exists an element 
$e\in A$ such that $\lt_e^\circ=\rt_e^\circ=\I_A$.
Now $e\circ x= \a(e)\b(x)$, that is, $\lt_e^\circ=\lt_{\a(e)}\beta$, so $\lt_e^\circ=\I_A$
if and only if $\lt_{\a(e)}$ is invertible $\b=\lt_{\a(e)}\inv$. Now $\lt_{\a(e)}$ being
invertible means that $\a(e)\in A^\ast$, and thus $\b=\lt_{\a(e)}\inv= \lt_{\a(e)\inv}$.
Similarly, $\rt_e^\circ=\rt_{\b(e)}\a$ equals the identity map if and only if
$\a=\rt_{\b(e)}\inv=\rt_{\b(e)\inv}$.

It readily verified that $(ab)\inv\circ x= x = x\circ (ab)$ in $A_{\rt_a,\lt_b}$.
\end{proof}

\begin{prop} \label{isomorphic}
  Let $A$ be a Hurwitz algebra. Any isotope $B$ of $A$ that has an identity element is
  again a Hurwitz  algebra, isomorphic to $A$, and $n_B=n_A(1_B)\inv n_A$.
\end{prop}

\begin{proof}
By Lemma~\ref{identity}, any principal isotope $B=(A,\circ)$ of $A$ that is unital has
the form $B=A_{\rt_c,\lt_d}$. Defining $n_B(x)=n_A(cd)n_A(x)$ for all $x\in B$,
we have
\begin{align*}
n_B(x\circ y)&= n_B((xc)(dy))= n_A(cd)n_A((xc)(dy)) \\
&= n_A(cd)n_A(x)n_A(cd)n_A(y)=n_B(x)n_B(y)
\end{align*}
so $B$ is a Hurwitz algebra. 
Moreover, $\lt_{cd}:(B,n_B)\to (A,n_A)$ is an isometry. Being isometric as quadratic
spaces, $A$ and $B$ are isomorphic algebras.
Since $1_B=(cd)\inv$, it is clear that $n_B=n_A(cd)n_A=n_A(1_B)\inv n_A$.
\end{proof}

From Proposition~\ref{isomorphic} follows the simple but important observation that, if
two isotopes $A_{\a,\b}$ and $B_{\gamma,\d}$ of Hurwitz algebras $A$ and $B$ are
isomorphic to each other, then $A\simeq B$.

\begin{cor} \label{similitud}
  Let $A$ be a Hurwitz algebra, and $\a,\b,\gamma,\d\in\gl(A)$.
  Any isomorphism $\f:A_{\a,\b}\to A_{\gamma,\d}$ is a similitude of $(A,n_A)$ with
  multiplier $n_A(\f(1_A))$.
\end{cor}

\begin{proof}
Let $\f:A_{\a,\b}\to A_{\gamma,\delta}$ be an isomorphism. It is straightforward to
verify that $\f$ is also an isomorphism $A\to B$, where 
$B= A_{\gamma\f\a\inv\f,\delta\f\b\inv\f}$ Thus, in particular, $\f$ is an
orthogonal map $(A,n_A)\to(B,n_B)$, and $\f(1_A)=1_B$.
By Proposition~\ref{isomorphic}, $n_B=n_A(1_B)\inv n_A$, so
$n_A(\f(x))=n_A(1_B)n_B(\f(x))=n_A(\f(1_A))n_A(x)$.
\end{proof}

The following proposition describes all isomorphisms between isotopes of Hurwitz algebras. 
It generalises the description of eight-dimensional absolute-valued algebras given in
\cite[Theorem~4.3]{ckmmrr11}.

\begin{prop} \label{isocriterion}
  Let $A=(A,n)$ be a Hurwitz algebra and $\a,\b,\gamma,\d\in\gl(A)$.
  A map $\f:A_{\a,\b}\to A_{\gamma,\d}$ is an isomorphism if and only if $\f\in \g(A)$ and
  there exist triality componenets $\f_1$, $\f_2$ of $\f$ such that 
    \begin{equation*}
      \begin{cases}
        \gamma= \f_1\a\f\inv, \\ 
        \delta= \f_2\b\f\inv.    
      \end{cases}
    \end{equation*}
\end{prop}

\begin{rmk}
  Note that, since any pair $(\f_1,\f_2)$ of triality components of $\f$ satisfy
  $\f_1=\rt_{\f_2(1)}\inv$, $\f_2=\lt_{\f_1(1)}\inv$, we could also have written
  \begin{equation*}
    \begin{cases}
      \gamma= \f_1\a\f\inv =  \rt_{\f_2(1_A)}\inv\f\a\f\inv, \\
      \delta= \f_2\b\f\inv =  \lt_{\f_1(1_A)}\inv\f\b\f\inv
    \end{cases}
  \end{equation*}
  in the proposition above.
\end{rmk}

\begin{proof}
A direct computation shows that if $\f\in\g(A)$ and $\gamma= \f_1\a\f\inv$ and
$\delta=\f_2\b\f\inv$, then $\f\in\Iso(A_{\a,\b}, A_{\gamma,\delta})$.

For the other direction, assume that $\f\in\Iso(A_{\a,\b},A_{\gamma,\d})$. This means that
$\f(\a(x)\b(y))=\gamma\f(x)\cdot\d\f(y)$ for all $x,y\in A$, and setting $z=\a(x)$,
$w=\b(y)$ gives 
$$\f(zw)=\gamma\f\a\inv(z)\cdot\d\f\b\inv(w).$$
Inserting $w=1_A$ into this equation yields 
$\f(z)=\gamma\f\a\inv(z)\cdot\d\f\b\inv(1_A)$, so 
$\gamma\f\a\inv= \rt_{\d\f\b\inv(1_A)}\inv\f$. By Corollary~\ref{similitud}, we know that
$\f$ is a similitude with multiplier $n_A(\f(1_A))$, so it follows that
$\gamma\f\a\inv$ is a similitude with multiplier
$\frac{n_A(\f(1_A))}{n_A(\d\f\b\inv(1_A))}$.
Similarly, $\d\f\b\inv$ is a similitude with multiplier
$\frac{n_A(\f(1_A))}{n_A(\gamma\f\a\inv(1_A))}$. 
Setting $\f_1=\gamma\f\a\inv$ and $\f_2=\d\f\b\inv$ gives $\f(zw)=\f_1(z)\f_2(w)$ with
$\f_1,\f_2\in\go(A)$, which means that $\f\in \g(A)$. Moreover, $\gamma= \f_1\a\f\inv$ and
$\delta=\f_2\b\f\inv$, as required.
\end{proof}

We record the following observations for future use. The first statement is a consequence
of the fact, referred to in the introduction, that $x\in N(A)$ if and only if $(ax)b=a(xb)$
for all $a,b\in A$; the second follows from Proposition~\ref{isocriterion}.

\begin{lma} \label{equalisotopes}
  Let $A$ be a Hurwitz algebra, $\a,\b,\gamma,\delta\in\gl(A)$ and $\rho\in k^\ast$.
  \begin{enumerate}
  \item $A_{\a,\b}=A_{\gamma,\delta}$ if and only if $\a=\rt_w\inv\gamma$,
    $\b=\lt_w\delta$ for some $w\in N(A)^\ast$.
  \item The homothety $h_\rho(x)= \rho x$ on $A$ defines an isomorphism $A_{\rho\I,\I}\to A$.
  \end{enumerate}
\end{lma}

The set 
$T=\{([\f_0],[\f_1],[\f_2])\in\pg(A)^3\mid \f_0(xy)=\f_1(x)\f_2(y)\,,\,\forall x,y\in A\}$
is a group under component-wise multiplication. The kernel of the group epimorphism
$$\pi_0:T\to\pg(A), (\f_0,\f_1,\f_2)\mapsto\f_0$$ 
is 
$$N=\{([\I_A],[\rt_w\inv],[\lt_w]) \mid w\in N(A)^\ast\}\simeq N(A)^\ast/(k^\ast1)\,,$$
hence $\pi_0$ induces an isomorphism $T\slash N\tilde{\to}\pg(A)$.
Now $T$ acts on $\pgl(A)^2$ by
\begin{equation} \label{preaction}
([\f_0],[\f_1],[\f_2])\cdot([\a],[\b])=([\f_1\a\f_0\inv],[\f_2\b\f_0\inv])\,.
\end{equation}
We denote by $X_A=\pgl(A)^2/N$ the orbit set of the induced $N$-action on
$\pgl(A)^2$. 
Then (\ref{preaction}) induces an action of $\pg(A)\simeq T\slash N$ on $X_A$, given by
\begin{equation} \label{groupaction}
  \f\cdot\left[\a,\b\right] = \left[\f_1\a\f\inv, \f_2\b\f\inv \right]
\end{equation}
for $[\a,\b]\in X_A$ and $\f\in\pg(A)$. Here $\f_1$, $\f_2$ are any pair of triality
components of $\f$, that is, $(\f,\f_1,\f_2)\in\pi_0\inv(\f)\subset T$.

We use the following notational convention:  
if $\gamma,\delta\in\gl(A)$, then the element of $X_A$ represented by
$([\gamma],[\delta])\in\pgl(A)$ is denoted by $[\gamma,\delta]$ 
(rather than $\left[[\gamma],[\delta]\right]$). 

Let $\mathscr{X}(A) = {}_{\pg(A)} \: X_A$ be the
groupoid of the action~(\ref{groupaction}).
For any Hurwitz algebra $A$, let $\mathscr{I}(A)$ denote the category of principal
isotopes of $A$, and $\check{\mathscr{I}}(A)$ the groupoid obtained from $\mathscr{I}(A)$
by removing all non-isomorphisms between the objects.
Note that if $A$ is a division algebra then so are all its isotopes, and thus any
non-zero morphism in $\mathscr{I}(A)$ is an isomorphism in this case.

The essence of our findings so far is summarised in the following theorem.

\begin{thm} \label{huvudsats}
  For any Hurwitz algebra $A$, the categories $\check{\mathscr{I}}(A)$ and
  $\mathscr{X}(A)$ are equivalent. 
  An equivalence $\mathscr{F}_{\!A}: \check{\mathscr{I}}(A) \to \mathscr{X}(A)$
  is given by $\mathscr{F}_{\!A}(A_{\a,\b})=[\a,\b]$ and $\mathscr{F}_{\!A}(\f)=[\f]$.
\end{thm}

\begin{proof}
Let $(\a,\b),(\gamma,\delta)\in\gl(A)^2$. If $A_{\a,\b}=A_{\gamma,\delta}$ then, by
Lemma~\ref{equalisotopes},
$$(\gamma,\delta)=(\rt_w\inv\a,\lt_w\b)=(\I_A,\rt_w\inv,\lt_w)\cdot(\a,\b)$$
for some  $w\in N(A)^\ast$, and hence the elements $([\a],[\b])$ and $([\gamma],[\d])$ in
$\pgl(A)^2$ represent the same object in $X_A=\pgl(A)^2\slash N$.
Proposition~\ref{isocriterion} guarantees that $[\f]\cdot[\a,\b]=[\gamma,\d]$ whenever
$\f\in\Iso(A_{\a,\b},A_{\gamma,\d})$. This shows that $\mathscr{F}_A$ is well defined. 
Clearly, $\mathscr{F}_A$ is surjective on objects, hence dense as a functor.

If $\f,\psi\in\Iso(A_{\a,\b},A_{\gamma,\delta})$ and
$\mathscr{F}_A(\f)=\mathscr{F}_A(\psi)$, then $\f=\rho\psi$ for some $\rho\in k^\ast$, so
$\rho\I_A=\f\psi\inv\in\Aut(A_{\a,\b})$. This implies $\rho=1$ and $\f=\psi$\,; hence
$\mathscr{F}_A$ is faithful.
Fullness is clear from the construction. 
\end{proof}

\begin{rmk} \label{anm}
\begin{enumerate}
\item
  If $k$ is a Euclidean field\footnote{A field $k$ is called Euclidean if $k^{\ast2}$
    forms an ordering of $k$.} 
  (\eg, $k=\R$), then the group $\pgo(A,n)$ is canonically isomorphic to
  $\og(A,n)/(\pm\I)$, via composition of inclusion and quotient projection:
  $\og(A,n)/(\pm\I)\hookrightarrow\go(A,n)/(\pm\I)\twoheadrightarrow\go(A,n)/(k^\ast\I)=\pgo(A,n)$.
  This also induces an isomorphism $\og^+(A,n)/(\pm\I)\to\pgo^+(A,n)$.
\item
  The set $X_A=\pgl(A)^2/N$ may also be viewed as the orbit set 
  $\pgl(A)^2/N(A)^\ast$ of the $N(A)^\ast$-action on $\pgl(A)^2$ given by
  $w\cdot\left([\a],[\b]\right)=\left([\rt_w\inv\a],[\lt_w\b] \right)$.
\item 
  Assume that $A$ is associative. The maps 
  $\pi_i:T\to\pg(A),\:(\f_0,\f_1,\f_2)\mapsto \f_i$, $i=1,2$, are group epimorphisms, and 
  since $N$ is normal in $T$, the images 
  $\pi_1(N)=\{[\rt_a]\mid a\in A^\ast\}\subset\pg(A)$ and 
  $\pi_2(N)=\{[\lt_a]\mid a\in A^\ast\}\subset\pg(A)$ of $N$ under $\pi_1$ and $\pi_2$ are
  normal subgroups of $\pg(A)$. It follows that $\rt_{A^\ast}=\{\rt_a\mid a\in A^\ast\}$ and
  $\lt_{A^\ast}=\{\lt_a\mid a\in A^\ast\}$ are normal subgroups of $\g(A)$.
\item
  If $\dim A=8$, since $N(A)=k1$, we have $X_A=\pgl(A)^2$. 
\end{enumerate}
\end{rmk}

If $A$ is associative, then the group $\g(A)$ has a particularly nice form.

\begin{prop} \label{asshurw}
  If $A$ is an associative Hurwitz algebra, then $\g(A)$ can be written as a semi-direct
  product in the following ways:
  $$\g(A)= \lt_{A^\ast}\rtimes\Aut(A)=\rt_{A^\ast}\rtimes\Aut(A) \,.$$
In particular, $\lt_{A^\ast}$ and $\rt_{A^\ast}$ are normal subgroups of $\g(A)$. 
\end{prop}

\begin{proof}
  As already noted in Remark~\ref{anm}(2), $\lt_{A^\ast}$ and $\rt_{A^\ast}$ are normal
  subgroups of $\g(A)$. 
  Since any automorphism of $A$ must fix the identity element $1_A$ in $A$, we have
  $\lt_a\in\Aut(A)$ for $a\in A^\ast$ if and only if $a=1$. Thus $\lt_{A^\ast}\cap\Aut(A)=\{\I\}$. 
  The inclusion $\Aut(A)\subset\g(A)$ is obvious.
  If $\f\in\g(A)$, then
  $$\lt_{\f(1)}\inv\f(xy)=\lt_{\f(1)}\inv\left(\f(x)\f(1)\inv\f(y)\right)=
  \lt_{\f(1)}\inv\f(x)\cdot\lt_{\f(1)\inv}\f(y) =
  \lt_{\f(1)}\inv\f(x)\cdot\lt_{\f(1)}\inv\f(y),$$ 
  so $\lt_{\f(1)}\inv\f\in\Aut(A)$. This implies $\g(A)=\lt_{A^\ast}\Aut(A)$, which
  shows that $\g(A)= \lt_{A^\ast}\rtimes\Aut(A)$.
  The identity $\g(A)= \rt_{A^\ast}\rtimes\Aut(A)$ is proved similarly.
\end{proof}

Note that if $\dim A=2$ then $\Aut(A)=\{\I,\kappa\}=\ct$.
If $A$ is a quaternion algebra then it is central simple, and the Skolem-Noether Theorem
\cite[p.~222]{basicII} gives $\Aut(A)=\{\lt_a\rt_a\inv\mid a\in A^\ast\}$.
Hence every $\f\in\go^+(A)$ can be written as 
$\f=\lt_a\lt_b\rt_b\inv=\lt_{ab}\rt_{b\inv}$.
This leads to the following result.

\begin{cor} \label{go4}
  Every proper similitude of a quaterion algebra $A$ has the form $\lt_a\rt_b$ for some 
  $a,b\in A^\ast$. The kernel of the group epimorphism 
  $A^\ast\times(A^\ast)^{\rm op}\to\go^+(A),\:(a,b)\mapsto\lt_a\rt_b$ is
  $\{(\rho,\rho\inv)\mid \rho\in k^\ast\}\subset A^\ast\times (A^\ast)^{\rm op}$. 
\end{cor}

\begin{proof}
From the argument in the paragraph preceding Corollary~\ref{go4}, it is clear that every $\f\in\go^+(A)$ can be written as $\lt_a\rt_b$ for some $a,b\in A^*$.
For the second assertion, assume that $\lt_a\rt_b=\I_A$, \ie, that $axb=x$ for all 
$x\in A$. Setting $x=a\inv$ gives $b=a\inv$, and now $\lt_{a}\rt_{a\inv}=\I_A$ if and only
if $a\in k^*1$.
\end{proof}

If $A$ is an associative Hurwitz algebra, then the map 
$\lt_a\psi\in \lt_{A^\ast}\rtimes \Aut(A)=\g(A)$ has a pair of triality components 
$(\lt_a\psi,\psi)$.
Thus the action (\ref{groupaction}) of 
$[\lt_a\psi]\in\pg(A)$ on $[\a,\b]\in X_A$ can be written as
\begin{equation} \label{groupaction2}
[\lt_a\psi]\cdot[\a,\b]=\left[\lt_a\psi\a\psi\inv\lt_a\inv,\psi\b\psi\inv\lt_a\inv\right] \,.
\end{equation}
For $\dim A=2$, this description of the groupoid $\mathscr{X}(A)$ is equivalent to the
isomorphism criterion 1.12 in \cite{petersson00}, applied to isotopes of quadratic \'etale
algebras.

We conclude this section by introducing a numerical, easily computed isomorphism invariant.
For $l$ a positive integer, let 
$k^{\ast l}=\{\rho^l\in k^\ast\mid \rho\in k^\ast\}\subset k^\ast$.

\begin{prop} \label{gensign}
  Let $A$ be a Hurwitz algebra of dimension $2l\ge4$, 
  and $\a,\b\in\gl(A)$.
  Then the pair $([\det(\a)],[\det(\b)])\in \left(k^\ast/k^{\ast l}\right)^2$ is an isomorphism
  invariant for $A_{\a,\b}$.
\end{prop}

\begin{proof}
Assuming the existence of an isomorphism $\f:A_{\a,\b}\to A_{\gamma,\delta}$, we need to show
that $\det(\a)\det(\gamma)\inv\in k^{*l}$ and $\det(\b)\det(\delta)\inv\in k^{*l}$.
For any  $x\in A$, we have
$\f\lt_{\a(x)}\b=\lt_{\gamma\f(x)}\d\f$ and hence 
$\f\lt_{\a(x)}\b\f\inv=\lt_{\gamma\f(x)}\d$, so
\begin{equation}\label{gensigneq}
\det\left(\lt_{\a(x)}\right)\det(\b)= \det\left(\f\lt_{\a(x)}\b\f\inv\right) =
\det\left( \lt_{\gamma\f(x)}\d\right) = \det(\lt_{\gamma\f(x)})\det(\d).
\end{equation}
Since the maps $\lt_a$ and $\rt_a$ are proper similitudes with multiplier
$n_A(a)$ for any $a\in A^\ast$, the equation \eqref{gensigneq} implies that 
$n_A(\a(x))^{l}\det(\b)= n_A(\gamma\f(x))^{l}\det(\d)$ whenever 
$\a(x)\in A^\ast$. Hence, $\det(\b)\det(\d)\inv\in k^{\ast l}$.
Similarly, the identity $\f\rt_{\b(y)}\a\f\inv=\rt_{\d\f(y)}\gamma$ implies that
$\det(\a)\det(\gamma)\inv\in k^{\ast l}$.
\end{proof}

Let $A$ and $l$ be as in Proposition~\ref{gensign}. For $i,j\in k^\ast/k^{\ast l}$, setting 
$$\mathscr{X}(A)_{i,j}= 
\left\{[\a,\b]\in\mathscr{X}(A)\mid\left([\det(\a)],[\det(\b)]\right)=(i,j)
\in 
\left(k^\ast/k^{\ast l}\right)^2\right\} \subset\mathscr{X}(A) \,,$$ 
the groupoid $\mathscr{X}(A)$ can be written as a coproduct 
\begin{equation} \label{gensigndecomp}
\mathscr{X}(A)=\coprod_{i,j\in\, k^\ast\!/k^{\ast l}}\mathscr{X}(A)_{i,j} \,.
\end{equation}
Hence, any subcategory $\mathscr{A}\subset\mathscr{X}(A)$ can be classified by classifying
each of the subcategories
$\mathscr{A}_{i,j}=\mathscr{A}\cap\mathscr{X}(A)_{i,j}\subset\mathscr{X}(A)_{i,j}$. 

As for the real ground field, $[\R^\ast:\R^{\ast l}]=2$ for any even integer $l\ge2$, the two
cosets being represented by $1$ and $-1$, and the quotient projection
$\R^\ast\to\R^\ast/\R^{\ast l}$ is given by $\rho\mapsto[\sign(\rho)]$.
If $A$ is either $\H$ or $\O$,
$\a,\b\in\gl(A)$ and $A_{\a,\b}=(A,\circ)$, then
$\det(\lt_x^\circ)=\det(\lt_{\a(x)})\det(\b)$. 
Since $\det(\lt_{\a(x)})=n(\a(x))^{(\dim A)/2}>0$ for any $x\ne0$, it follows that
$\sign(\det(\lt^\circ_x))=\sign(\det(\b))$, and similarly
$\sign(\det(\rt_x^\circ))=\sign(\det(\a))$. 
This means that the decomposition
$\mathscr{X}(A)=\coprod_{i,j\in\{-1,1\}}\mathscr{X}(A)_{i,j}$ here coincides with
the ``double sign'' decomposition for real division algebras, introduced in
\cite{doublesign}.

\section{Isotopes of Hamilton's quaternions} \label{euclidean}

In this section , we give a more detailed account for isotopes of
real Hurwitz algebras whose underlying quadratic space is Euclidean.
Such a Hurwitz algebra is isomorphic to either $\R$, $\C$, $\H$, or $\O$, and the isotopes
are precisely the real division algebras that are isotopic to a Hurwitz algebra.
The isotopes of $\C$ are all the two-dimensional real division algebras, and their
classification has been described in \cite{dieterich05,hupe04}.
While our goal is to study the isotopes of $\H$, some arguments hold also in the
$8$-dimensional case. Thus, let $A$ be either $\H$ or $\O$. 
Our main tool will be polar decomposition of linear maps.

Any $\a\in\gl(A)$ can be written
as $\a=\a'\l$, where $\det(\a')=|\det(\a)|>0$ and $\l\in \ct=\{\I,\kappa\}$. Polar
decomposition now yields $\a'=\zeta\d$, with $\zeta\in\so(A)$ and $\d\in\Pds(A)$.
Hence $\a=\zeta\d\l$, and this decomposition is unique \cite[\S 14]{gantmacher59i}.

Passing to the projective setting, the above implies that every $\a\in\pgl(A)$ uniquely
determines $\zeta\in\so(A)/(\pm\I)$, $\d\in\spds(A)=\Pds(A)\cap\SL(A)$ and $\l\in \ct$
such that $\a=[\zeta\d\l]$.
As noted in Remark~\ref{anm}, $\so(A)/(\pm\I)\simeq\pgo^+(A)$, and $\zeta$ can indeed be
viewed as an element in $\pgo^+(A)$ instead of $\so(A)/(\pm\I)$.
The cyclic group $\ct$ acts on $\pgo^+(A)$ by $[\f]^\l=[\l\f\l]$
($\l\in\ct$, $\f\in\go^+(A)$). 
Note that $[\lt_a]^\kappa=[\rt_{\bar{a}}]=[\rt_a\inv]$ and
$[\rt_a]^\kappa=[\lt_{\bar{a}}]=[\lt_a\inv]$ in $\pgo^+(A)$.
We write $[\f]^{-\l}$ for $\left([\f]\inv\right)^\l=\left([\f]^\l\right)\inv$.

Let $\a,\b\in\pgl(A)$, $\a=[\zeta\d\l]$ and $\b=[\eta\e\mu]$, where $\zeta,\eta\in\pgo^+(A)$,
$\d,\e\in\spds(A)$ and $\l,\mu\in \ct$.
The action (\ref{groupaction}) of $[\f]\in\pgo^+(A)$ on $[\a,\b]\in X_A$ is given by
\begin{equation*}
[\f]\cdot[\a,\b]=\left[\f_1\zeta\d\l\f\inv \,,\, \f_2\eta\e\mu\f\inv\right] =
\left[\f_1\zeta\f^{-\l}(\f^\l\d\f^{-\l})\l \,,\, \f_2\eta\f^{-\mu}(\f^\mu\e\f^{-\mu})\mu\right] \,,
\end{equation*}
and $\f_1\zeta\f^{-\l},\,\f_2\eta\f^{-\mu}\in\pgo^+(A)$\,,
$\f^\l\d\f^{-\l},\,\f^\mu\e\f^{-\mu}\in\spds(A)$.

The group $N(A)^{\ast}$ acts on $\pgo^+(A)^2$ by 
$w\cdot(\zeta,\eta)=([\rt_w]\inv\zeta\,,\,[\lt_w]\eta)$,
and we denote the orbit set of this action by $\pgo^+(A)^2/N(A)^{\ast}$.
The argument in the preceding paragraph shows that the $\pgo^+(A)$-action on
$\pgo^+(A)^2/N(A)^\ast \times \spds(A)^2 \times \ct^2$ by
\begin{equation} \label{factorisedaction}
[\f]\cdot\left( [\zeta,\eta], (\d,\e), (\l,\mu) \right) = 
\left( \left[\f_1\zeta\f^{-\l},\f_2\eta\f^{-\mu}\right], (\f^\l\d\f^{-\l},\f^\mu\e\f^{-\mu}),
  (\l,\mu)\right)
\end{equation}
is equivalent to the $\pgo^+(A)$-action (\ref{groupaction}) on $X_A$, via the map
$$m:\pgo^+(A)^2/N(A)^\ast \times \spds(A)^2 \times \ct^2\to X_A \,,\;
\left( [\zeta,\eta], (\d,\e), (\l,\mu) \right)\mapsto [\zeta\d\l, \eta\e\mu]$$
(\ie, $\f m=m\f$ for all $\f\in\pgo^+(A)$). Hence, in particular, the groupoid 
$$\mathscr{Y}(A)=
{}_{\pgo^+(A)} \left(\: \pgo^+(A)^2/N(A)^\ast \times \spds(A)^2 \times \ct^2 \:\right)$$
of the action (\ref{factorisedaction}) is isomorphic to $\mathscr{X}(A)$:

\begin{prop}
  The functor
  $ \mathscr{G}_A:\mathscr{Y}(A) \to \mathscr{X}(A)$
  defined by $\mathscr{G}_A(y)=m(y)$ for $y\in\mathscr{Y}(A)$, and $\mathscr{G}_A(\f)=\f$
  for morphisms $\f\in\pgo^+(A)$, is an isomorphism of categories.
\end{prop}

From here on we focus exclusively on the four-dimensional case, in which $A\simeq\H$.
Corollary~\ref{go4} implies that every element $\zeta\in\pgo^+(\H)$ has the form
$\zeta=[\lt_a\rt_b]$ for some $a,b\in\H^\ast$. Since $\H$ is associative,
it follows that if $\f=[\lt_a\rt_b]$ then $\f_1=[\lt_a]$, $\f_2=[\rt_b]$ are triality
components of $\f$. 
Recall that $N(\H)^\ast=\H^\ast$ acts on $\pgo^+(\H)^2$ by
$c\cdot(\zeta,\eta)=\left([\rt_c\inv]\zeta,[\lt_c]\eta\right)$.
Since $\lt_a\rt_b=\rt_b\lt_a$ for all $a,b\in\H$, the following identities hold
in $\pgo^+(\H)^2/N(\H)^\ast=\pgo^+(\H)^2/\H^\ast$:
$$[\lt_a\rt_b,\lt_c\rt_d]=c\inv\cdot[\lt_a\rt_b,\lt_c\rt_d]=
[\rt_c\lt_a\rt_b,\rt_d]=[\lt_a\rt_{bc},\rt_d]\,.$$
Hence every element
$[\zeta,\eta]\in\pgo^+(\H)^2/\H^\ast$ can be written on the form
$[\zeta,\eta]=[\lt_a\rt_b,\rt_c]$, with $a,b,c\in\H^\ast$.
Now the action of $\f=[\lt_s\rt_t]\in\pgo^+(\H)$ ($s,t\in\H^\ast$) on
$([\zeta,\eta], (\d,\e), (\l,\mu))\in\mathscr{Y}(\H)$ is given by
\begin{equation}
  \f\cdot ([\zeta,\eta], (\d,\e), (\l,\mu))
  =\left( 
  \left[\lt_s\lt_a\rt_b\lt_s^{-\l}\rt_t^{-\l},
    \rt_t\rt_c\lt_s^{-\mu}\rt_t^{-\mu}\right],\,
    \left(\f^\l\d\f^{-\l}, \f^\mu\e\f^{-\mu}\right),\, (\l,\mu) \right).
\end{equation}

For $(i,j)\in\{-1,1\}^2$, set 
$\mathscr{Y}(\H)_{i,j}=\mathscr{G}_A\inv\left(\mathscr{X}(\H)_{i,j}\right)
\subset\mathscr{Y}(\H)$.
Note that $([\zeta,\eta], (\d,\e),
(\kappa^i,\kappa^j))\in\mathscr{Y}(\H)_{(-1)^i,(-1)^j}$.
This gives the decomposition 
$$\mathscr{Y}(\H)=\coprod_{i,j\in\{-1,1\}}\mathscr{Y}(\H)_{i,j}\,,$$
and the action of $\f=\lt_s\rt_t$ on each of the cofactors $\mathscr{Y}(\H)_{i,j}$ can be
studied separately.

For $([\zeta,\eta], (\d,\e), (\I,\I))\in\mathscr{X}(\H)_{1,1}$, we have
\begin{align*}
  \f\cdot ([\zeta,\eta], (\d,\e), (\I,\I))
  &=\left(\left[\lt_s\lt_a\rt_b\lt_s\inv\rt_t\inv,\rt_t\rt_c\lt_s\inv\rt_t\inv\right]
    \,,\, \left(\f\d\f\inv, \f\e\f\inv\right) \,,\, (\I,\I) \right) \\
  &=\left(\left[\rt_s\inv\rt_b\rt_t\inv\lt_s\lt_a\lt_s\inv,\rt_t\rt_c\rt_t\inv\right]
    \,,\, \left(\f\d\f\inv, \f\e\f\inv\right) \,,\, (\I,\I) \right) \,.
\end{align*}
In particular, if $s=1$, $t=b$, so that $\f=\rt_b$, then
$$\f\cdot([\zeta,\eta], (\d,\e), (\I,\I))  =\left(\left[\lt_a,\rt_{bcb\inv}\right]
    \,,\, \left(\rt_b\d\rt_b\inv, \rt_b\e\rt_b\inv\right) \,,\, (\I,\I) \right) \,.$$
Hence every orbit in $\mathscr{Y}(\H)_{1,1}$ contains an element of the form
$([\lt_a,\rt_b], (\d,\e), (\I,\I))$. Moreover, the action of $\f=[\lt_s\rt_t]$ on such an
element is given by
\begin{equation} \label{stabilising}
\begin{split}
  \f\cdot([\lt_a,\rt_b], (\d,\e), (\I,\I))&=
\left(\left[\lt_s\lt_a\lt_s\inv\rt_t\inv, \rt_t\rt_b\lt_s\inv\rt_t\inv\right],
  (\f\d\f\inv,\f\e\f\inv), (\I,\I)\right) \\
&=\left(\left[\rt_{(st)\inv}\lt_{sas\inv},\rt_{t\inv bt}\right],
  (\f\d\f\inv,\f\e\f\inv), (\I,\I)\right)
\end{split}
\end{equation}
which has the form $([\lt_c,\rt_d], (\d,\e), (\I,\I))$ if and only if $t=s\inv$, that
is, $\f=[\lt_s\rt_s\inv]$.

Denoting by $\mathscr{Z}={}_{\H^\ast/\R^\ast}\left((\H^\ast/\R^\ast)^2\times
\spds(\H)^2\right)$ the groupoid of the action
$$s\cdot((a,b),(\d,\e)) = ((sas\inv,sbs\inv), (c_s\d c_s\inv,c_s\e c_s\inv)) \,,$$
the above implies that the functor 
$\mathscr{H}_{1,1}:\mathscr{Z}\to \mathscr{Y}(\H)_{1,1}$ defined by 
\begin{align*}
  &\mathscr{H}_{1,1}((a,b),(\d,\e)) = ([\lt_a,\rt_b], (\d,\e),(\I,\I))\\
\intertext{on objects, and}
  &\mathscr{H}_{i,j}(s)=[\lt_s\rt_s\inv]
\end{align*}
on morphisms, is an equivalence of categories.

Next, consider $([\zeta,\eta],(\d,\e),(\kappa,\I))\in\mathscr{Y}(\H)_{-1,1}$. In this case,
\begin{align*}
  \f\cdot([\zeta,\eta],(\d,\e),(\kappa,\I))
  &=\left(\left[\lt_s\lt_a\rt_b\lt_s^{-\kappa}\rt_t^{-\kappa},\rt_t\rt_c\lt_s\inv\rt_t\inv\right]
    \,,\, \left(\f^\kappa\d\f^{-\kappa}, \f\e\f\inv\right) \,,\, (\kappa,\I) \right) \\
  &=\left(\left[\lt_s\lt_a\lt_t\rt_s\inv\rt_b\rt_s ,\rt_t\rt_c\rt_t\inv\right]
    \,,\, \left(\f^\kappa\d\f^{-\kappa}, \f\e\f\inv\right) \,,\, (\kappa,\I) \right)
\end{align*}
Setting $s=1$ and $t=a\inv$ gives $\f=\rt_a\inv$ and
$$[\rt_a]\inv\cdot([\zeta,\eta],(\d,\e),(\kappa,\I))=
\left(\left[\rt_b,\rt_{aca\inv}\right]
     \,,\, \left(\lt_a\d\lt_a\inv, \rt_a\inv\e\rt_a\right) \,,\, (\kappa,\I) \right)$$
so the orbit of $([\zeta,\eta],(\d,\e),(\kappa,\I))$ contains an element of the form
$\left([\rt_a,\rt_b],(\d',\e'),(\kappa,\I) \right)$.
Similarly to the previous case with $\mathscr{Y}(\H)_{1,1}$, the group elements
$\f\in\pgo^+(\H)$ stabilising this form, so that 
$\f\cdot([\rt_a,\rt_b],(\d,\e),(\kappa,\I))=([\rt_c,\rt_d],(\d',\e'),(\kappa,\I))$ for
some $c,d\in\H$, are precisely those of the form $\f=[\lt_s\rt_s\inv]$, $s\in\H^\ast$.
Note that if $\f=[\lt_s\rt_s\inv]$ then $\f^\kappa=\f$, hence 
$$\f\cdot([\rt_a,\rt_b],(\d,\e),(\kappa,\I))=
\left([\rt_{sas\inv},\rt_{sbs\inv}],(\f\d\f\inv,\f\e\f\inv),(\kappa,\I)\right) \,.$$
Again, this shows that there is an equivalence of categories 
$\mathscr{H}_{-1,1}:\mathscr{Z} \to \mathscr{Y}(\H)_{-1,1}$ is given by 
$\mathscr{H}_{-1,1}(s)=[\lt_s\rt_s\inv]$ on morphisms and
$\mathscr{H}_{-1,1}((a,b),(\d,\e)) = ([\rt_a,\rt_b], (\d,\e),(\kappa,\I))$
on objects.

Similar computations can be made for $\mathscr{Y}(\H)_{1,-1}$ and
$\mathscr{Y}(\H)_{-1,-1}$, in each case resulting in an equivalence with the category
$\mathscr{Z}$. The results are summarised in Theorem~\ref{quaternion} below.

For $a\in\H^\ast$, set $c_a=[\lt_a\rt_a\inv]\in\pgo^+(\H)$. Note that the kernel of the
morphism $\H^\ast\to\pgo^+(\H),\; a\mapsto c_a$ is $\R^\ast1$.

\begin{thm}\label{quaternion}
  Let $(i,j)\in \{-1,1\}^2$. Each of the subcategories $\mathscr{Y}(\H)_{i,j}$ of
  $\mathscr{Y}(\H)$ is equivalent to the groupoid 
  $\mathscr{Z}={}_{\H^\ast/\R^\ast}\left((\H^\ast/\R^\ast)^2\times \spds(\H)^2\right)$ of the action 
  $$s\cdot((a,b),(\d,\e)) =
  ((sas\inv,sbs\inv), (c_s\d c_s\inv,c_s\e c_s\inv)) \,.$$
  An equivalence 
$\mathscr{H}_{i,j}:\mathscr{Z} \to \mathscr{Y}(\H)_{i,j}$ is given by 
$\mathscr{H}_{i,j}(s)=c_s$ for morphisms
$s\in\H^\ast/\R^\ast$ and
\begin{align*}
  \mathscr{H}_{1,1}((a,b),(\d,\e)) &= ([\lt_a,\rt_b], (\d,\e),(\I,\I)) \,, &
  \mathscr{H}_{-1,1}((a,b),(\d,\e)) &= ([\rt_a,\rt_b], (\d,\e),(\kappa,\I)) \,, \\
  \mathscr{H}_{1,-1}((a,b),(\d,\e)) &= ([\lt_a,\lt_b], (\d,\e),(\I,\kappa)) \,, &
  \mathscr{H}_{-1,-1}((a,b),(\d,\e)) &= ([\lt_a,\rt_b], (\d,\e),(\kappa,\kappa)) \,.
\end{align*}
\end{thm}

\begin{rmk}
  Theorem~\ref{quaternion} shows, in particular, that
  $\mathscr{X}(\H)_{i,j}\simeq\mathscr{X}(\H)_{i',j'}$ for all
  $(i,j),(i',j')\in\{-1,1\}^2$. 
\end{rmk}

The image of the monomorphism $\H^\ast/\R^\ast\to\pgo^+(\H),\;s\mapsto c_s$ is
$\{\f\in\pgo^+(\H)\mid \f(\R^\ast1)=\R^\ast1\}$, which can be identified with
$\so_1(\H)=\{\f\in\so(\H)\mid \f(1)=1\}\simeq\so(1_{\H}^\perp)$. Thus the map
$f:\H^\ast/\R^\ast\to\so_1(\H),\;s\mapsto c_s$ is an isomorphism, inducing an action of
$\so_1(\H)$ on the object set $(\H^\ast/\R^\ast)^2\times \spds(\H)^2$ of $\mathscr{Z}$,
given by 
$\f\cdot((a,b),(\d,\e))=
f\inv(\f)\cdot((a,b),(\d,\e))=((\f(a),\f(b)),(\f\d\f\inv,\f\e\f\inv))$.

This allows for the following geometric interpretation of the category $\mathscr{Z}$.
Elements in $\H^\ast/\R^\ast$ are viewed as lines through the origin in $\H$ (that is,
elements in the real projective space $\mathbb{P}(\H)$), and
$\d\in\spds(A)$ is identified with the three-dimensional hyper-ellipsoid
$E_\d=\{x\in\H \mid \<x,\d(x)\>=1\}\subset\H$. The set
$\mathcal{E}=\{E_\d \mid \d\in\spds(\H)\}$ consists of all hyper-ellipsoids
centered in the origin and satisfying $\l_1\l_2\l_3\l_4=1$, where $\l_i\in\R_{>0}$,
$i=1,2,3,4$ are the lengths of the principal axes (which lie along the eigenspaces of $\d$).
Moreover, $\H$ is identified with
$\R^4$ in the natural way, and $1_\H^\perp\subset\H$ with $V=\spann\{e_2,e_3,e_4\}\subset\R^4$.
Now objects in $\mathscr{Z}$ can be viewed as configurations in $\R^4$ consisting of two
lines $a,b\in\mathbb{P}(\R^4)$, and two hyper-ellipsoids $E_\d,E_\e\in\mathcal{E}$.
A morphism $(a,b,E_\d,E_\e)\to(a',b',E_{\d'},E_{\e'})$ between two such configurations is
an element $\f\in\so(V)\subset\so(\R^4)$ transforming one configuration to the other:
$(\f(a),\f(b),\f(E_\d),\f(E_\e))=(a',b',E_{\d'},E_{\e'})$.

The above interpretation gives a picture of the category of isotopes of $\H$ in elementary
geometric terms, making it easy to read of whether or not two different isotopes of $\H$
are isomorphic. To extract from this description a complete and irredundant list of
representatives for the set of isomorphism classes of isotopes of $\H$, one would need to
find a normal form for the described configurations of lines and hyper-ellipsoids under
the natural action of $\so(V)$.
While this task is in principle not difficult, it is by its nature quite technical, and we
shall not pursue it further here. 
We note, however, that the object set of the groupoid $\mathscr{Z}$ is a $24$-dimensional
manifold, with isomorphism classes being orbits under a continuous action of a
$3$-dimensional real Lie group. 

\section{Composition algebras} \label{calg}

In this section, as an example, we show how the theory of composition algebras fits into
the general setting of Section~\ref{general}. This leads to a short and easy proof of
results by Stampfli-Rollier \cite{stampfli83} on the structure four-dimensional composition
algebras. As a bonus, our approach works also in characteristic two.

\begin{lma} \label{compgo}
  Let $A=(A,n_A)$ be a Hurwitz algebra, and $\a,\b\in\gl(A)$. The isotope $A_{\a,\b}$ is a
  composition algebra if and only if $\a,\b\in\go(A,n_A)$.
\end{lma}

\begin{proof}
If $\a,\b$ are similitudes, then $A_{\a,\b}$ is a composition algebra with respect to the
norm $n=\mu(\a)\mu(\b)n_A$.
Conversely, suppose that $A_{\a,\b}=(A,\circ)$ is a composition algebra with norm
$n$. Let $a\in A$ be an element such that $n(a)\ne0$, then $b=n(a)\inv(a\circ a)$ has norm
$1$, and thus the linear endomorphisms $\lt_b^{\circ}$ and $\rt_b^{\circ}$ of $A_{\a,\b}$
are orthogonal with respect to $n$. 
Now the isotope $B=\left(A_{\a,\b}\right)_{(\rt_b^{\circ})\inv,(\lt_b^{\circ})\inv}$ is a
Hurwitz algebra with norm $n$ and identity element $b\circ b$.\footnote{%
This argument, by which any finite-dimensional composition algebra is shown to be isotopic
to a Hurwitz algebra with the same norm, is sometimes referred to as \emph{Kaplansky's
  trick}, after \cite[p.~957]{kaplansky53}.}
It follows that $B$ is isotopic to $A$
and thus, by Proposition~\ref{isomorphic}, $n=\rho n_A$ for some $\rho\in k^\ast$. 
Hence $$n_A(\a(x))n_A(\b(y))=n_A(\a(x)\b(y))=n_A(x\circ y)=\rho\inv n(x\circ y)
=\rho\inv n(x)n(y)=\rho n_A(x)n_A(y)$$
for all $x,y\in A$. Inserting $y=1_A$ into this equation gives 
$n_A(\a(x))=\rho n_A(\b(1))\inv n_A(x)$, so $\a$ is a similitude with muliplier 
$\rho n_A(\b(1))\inv$ with respect to $n_A$. Similarly, $\b$ is a similitude with
multiplier $\rho n_A(\a(1))\inv$.
\end{proof}

Given a Hurwitz algebra $A$, let
$\mathscr{X}^c(A)\subset\mathscr{X}(A)$ be the full subcategory formed by all images
$[\a,\b]=\mathscr{F}_A(A_{\a,\b})$ where $A_{\a,\b}$ is a composition algebra. 
Lemma~\ref{compgo} implies that 
$\mathscr{X}^c(A)=\pgo(A,n_A)^2/N\subset\mathscr{X}(A)$.

For $\dim A\ge2$, the property of being a proper respectively improper similitude
is retained by the factors $\a,\b$ of $[\a,\b]\in\mathscr{X}^c(A)$ under the action of
$\pg(A)$.
This means that for each $(i,j)\in\{-1,1\}^2$, the subset
$$\mathscr{X}^{c(i,j)}(A)=\left(\pgo^i(A)\times\pgo^j(A)\right)/N(A)^\ast
\subset\mathscr{X}^c(A)$$
(where we use the notational convention $\pgo^{\pm1}(A)=\pgo^\pm(A)$)
is invariant under this action, and the category $\mathscr{X}^c(A)$ decomposes as
$$\mathscr{X}^c(A)=\coprod_{i,j\in\{-1,1\}}\mathscr{X}^{c(i,j)}(A)\,.$$
This refines the decomposition (\ref{gensigndecomp}) for $\dim A\ge4$:
$\mathscr{X}^{c(i,j)}(A)\subset\mathscr{X}^c(A)_{i,j}$ for all $i,j\in\{-1,1\}$. 
If $-1\not\in k^l$, $l=(\dim A)/2$, then $\mathscr{X}^{c(i,j)}(A)=\mathscr{X}^c(A)_{i,j}$,
but whereas the four sets $\mathscr{X}^{c(i,j)}(A)$ are always distinct,
$\mathscr{X}^c(A)_{i,j}=\mathscr{X}^c(A)$ for all $i,j\in\{-1,1\}$ in case $-1\in k^l$.

If $(A,n_A)$ is a Euclidean space, the groups $\pgo(A)$ and $\pgo^+(A)$ are canonically
identified with $\og(A)/(\pm\I_A)$ and $\so(A)/(\pm\I_A)$ respectively.
This, together with Theorem~\ref{huvudsats}, gives a description of all finite-dimensional
absolute valued algebras equivalent to Theorem~4.3 in \cite{ckmmrr11}.

We proceed to describe all composition algebras isotopic to a fixed quaternion algebra
$A=(A,n)$ over $k$. 
Every $[\a,\b]\in\mathscr{X}^c(A)$ can be written as 
$[\a,\b]=[\zeta\l,\eta\mu]$ with $\zeta,\eta\in\go^+(A)$ and $\l,\mu\in \ct$. 
Corollary~\ref{go4} now gives $\zeta=[\lt_a\rt_b]$, $\eta=[\lt_c\rt_d]$ for some
$a,b,c,d\in A^\ast$, and since $N(A)=A$, we have
$[\a,\b]=[\lt_a\rt_b\l,\lt_c\rt_d\mu]=[\lt_a\rt_{bc}\l,\rt_d]$ in $X_A$.
Thus, analogously with the Euclidean case treated in Section~\ref{euclidean}, $[\a,\b]$
can be written as $[\a,\b]=[\lt_a\rt_b\l,\rt_c\mu]$. Moreover, one reads off that
$[\a,\b]\in\mathscr{X}^{c(\det(\l),\det(\mu))}(A)$. 

Let $[\a,\b]=[\lt_a\rt_b,\rt_c]\in\mathscr{X}^{c(1,1)}(A)$. Now
$[\rt_b]\cdot[\a,\b]=[\lt_a,\rt_{b\inv cb}]$, so every $\pgo^+(A)$-orbit in
$\mathscr{X}^{c(1,1)}(A)$ contains an element of the form $[\lt_b,\rt_b]$.
Similarly, every orbit in $\mathscr{X}^{c(-1,1)}(A)$ contains an element of the form
$[\rt_a\kappa,\rt_b]$, every orbit in $\mathscr{X}^{c(1,-1)}(A)$ contains some
$[\lt_a,\lt_b\kappa]$, and every orbit in $\mathscr{X}^{c(-1,-1)}(A)$ contains an element of
the form $[\lt_a\kappa,\rt_b\kappa]$.
By computations similar to (\ref{stabilising}), one proves that $\f\in\pgo^+(A)$
stabilises each of these forms if and only if $\f=c_s$ for some $s\in A^\ast$, and that
\begin{align*}
  c_s\cdot[\lt_a,\rt_b]&=
  [\lt_s\lt_a\lt_s\inv,\rt_s\inv\rt_b\rt_s]=[\lt_{sas\inv},\rt_{sbs\inv}]\,,\\
  c_s\cdot[\rt_a\kappa,\rt_b]&=
  [\rt_s\inv\rt_a\rt_s\kappa,\rt_s\inv\rt_b\rt_s]=[\rt_{sas\inv}\kappa,\rt_{sbs\inv}]\,,\\
  c_s\cdot[\lt_a,\lt_b\kappa]&=
  [\lt_s\lt_a\lt_s\inv,\lt_s\lt_b\lt_s\inv\kappa]=[\lt_{sas\inv},\lt_{sbs\inv}\kappa]\,,\\
  c_s\cdot[\lt_a\kappa,\rt_b\kappa]&=
  [\lt_s\lt_a\lt_s\inv\kappa,\rt_s\inv\rt_b\rt_s\kappa]=[\lt_{sas\inv}\kappa,\rt_{sbs\inv}\kappa]\,.
  \end{align*}
This proves the following result, which gives an isomorphism criterion for all
composition algebras isotopic to $A$.

\begin{prop} \label{compositionalg}
  Let $A$ be a quaternion algebra. 
  Each of the categories $\mathscr{X}^{c(i,j)}(A)$ is equivalent to the groupoid
  $\mathscr{Z}^c(A)={}_{A^\ast/k^\ast}\left(A^\ast/k^\ast \right)^2$ of the action
  $s\cdot(a,b)=(sas\inv,sbs\inv)$.
  Equivalences $\mathscr{H}^c_{i,j}:\mathscr{Z}^c(A)\to\mathscr{X}^c(A)_{i,j}$ are given
  by $\mathscr{H}^c_{i,j}(s)=c_s$ for morphisms, and
  \begin{align*}
\mathscr{H}^c_{1,1}(a,b)&=[\lt_a,\rt_b] \,, &
  \mathscr{H}^c_{-1,1}(a,b)&=[\rt_a,\rt_b] \,, \\
  \mathscr{H}^c_{1,-1}(a,b)&=[\lt_a,\lt_b] \,, &
  \mathscr{H}^c_{-1,-1}(a,b)&=[\lt_a,\rt_b] \,.
\end{align*}
\end{prop}

A characterisation of all four-dimensional composition algebras over $k$ for $\chr k\ne2$,
similar to Proposition~\ref{compositionalg}, is given by Stampfli-Rollier in Section~3--4 of
\cite{stampfli83}. Her exposition also contains explicit isomorphism criteria in terms of
the parameters $a,b$.
The Euclidean case, comprising all composition algebras isotopic to $\H$, that is, 
all four-dimensional absolute valued algebras, has also been described by Ram\'irez
\'Alvarez \cite{ramirez99}. 
Forsberg, in his master thesis \cite{forsberg09}, refined that description to give an
explicit cross-section for the isomorphism classes of these algebras.

\bibliographystyle{plain}
\bibliography{../litt}

\end{document}